\documentclass[a4paper,12pt]{amsart}
\usepackage{amssymb}
\usepackage{ifthen}
\usepackage{graphicx}
\nonstopmode
\newtheorem{thm}{Theorem}
\newtheorem{cor}{Corollary}

\theoremstyle{definition}

\newenvironment{rem}{%
\bigskip
\noindent \textsl{{\sl Remark. }}}{\bigskip}

\newenvironment{pf}[1][]{%
 \vskip 3mm
 \noindent
 \ifthenelse{\equal{#1}{}}%
  {{\slshape Proof. }}%
  {{\slshape #1.} }%
 }%
{\qed\bigskip}

\newcounter{alphabet}
\newcounter{tmp}

\makeatletter
\newcommand{\Ref}[1]{\@ifundefined{r@#1}{}{\setcounter{tmp}{\ref{#1}}\Alph{tmp}}}
\makeatother
\newcounter{minutes}
\divide\time by 60
\newcounter{hours}
\multiply\time by 60
\addtocounter{minutes}{-\time}
\newcommand{\bpf}{\begin{pf}}
\newcommand{\epf}{\end{pf}}
\newcommand{\A}{{\mathcal A}}
\newcommand{\B}{{\mathcal B}}
\newcommand{\Bers}{{\mathcal R}}

\newcommand{\es}{{\mathcal S}}

\newcommand{\C}{{\mathbb C}}
\newcommand{\K}{{\mathcal K}}

\newcommand{\N}{{\mathcal N}}
\newcommand{\M}{{\mathcal M}}

\newcommand{\ID}{{\mathbb D}}
\newcommand{\D}{{\mathbb D}}

\newcommand{\inv}{^{-1}}

\begin{document}
%\hfill{File}

\bibliographystyle{amsplain}

%
%\begin{center}
%{\tiny \texttt{FILE:~\jobname .tex,
%        printed: \number\year-\number\month-\number\day,
%        \thehours.\ifnum\theminutes<10{0}\fi\theminutes}
%}
%\end{center}

\title[Radius problems associated with (pre-)Schwarzian derivative]%
{Radius problems associated with pre-Schwarzian and Schwarzian
derivatives
}

%%%%%%%% BEGIN TIMESTAMP
\def\thefootnote{}
\footnotetext{ \texttt{\tiny File:~\jobname .tex,
          printed: \number\day-\number\month-\number\year,
          \thehours.\ifnum\theminutes<10{0}\fi\theminutes}
} \makeatletter\def\thefootnote{\@arabic\c@footnote}\makeatother
%%%%%%%% END TIMESTAMP

\author{S. Ponnusamy}
\address{Department of Mathematics,
Indian Institute of Technology Madras,
Chennai - 600 036, India}
\email{samy@iitm.ac.in}

\author{S.K. Sahoo}
\address{Department of Mathematics,
Indian Institute of Technology Indore,
M-Block, I.E.T. DAVV Campus,
Khandwa Road, Indore 452 017, M.P., India}
\email{swadesh@iiti.ac.in}

\author{T. Sugawa${}^{~\mathbf{*}}$}
\address{Graduate School of Information Sciences,
Tohoku University, Aoba-ku, Sendai 980-8579, Japan}
\email{sugawa@math.is.tohoku.ac.jp}


\begin{abstract}
Some of important univalence criteria for a non-constant meromorphic function $f(z)$
on the unit disk $\ID$ involve its pre-Schwarzian or Schwarzian derivative.
We consider an appropriate norm for the pre-Schwarzian derivative, and
discuss the problem of finding the largest possible $r\in (0,1)$ for which
the pre-Schwarzian norm of the dilation $r^{-1}f(rz)$ is not greater than
a prescribed number for normalized univalent functions $f(z)$ in the unit disk.
Similar results concerning the Schwarzian derivative are also obtained.
%Some of the important univalence criteria for a non-constant meromorphic function $f(z)$
%on the unit disk $\ID$ involves its pre-Schwarzian and Schwarzian derivatives, that are defined by
%$T_f(z)$ and $S_f(z)$, respectively. For $r\in (0,1)$, consider the dilation $f_r(z)=r^{-1}f(rz)$,
%and for $\alpha\ge0$ and $M>0,$ we set
%$$\B_\alpha(M)=\{f\in \A: \|T_f\|_\alpha\le M\}
%~\mbox{ and }~\N_\alpha(M)=\{f\in \M: \|S_f\|_\alpha\le M\},
%$$
%where $\A$ (resp.~$\M$) denote the set of normalize analytic (resp.~meromorphic) functions $f$ on $\D$,
%and
%$$\|\varphi\|_\alpha=\sup_{z\in\D}(1-|z|^2)^\alpha|\varphi(z)|
%$$
%for a meromorphic function $\varphi$ on $\D.$ In this article we discuss certain radius properties of
%$\B_\alpha(M)$ and $\N_\alpha(M)$. As applications, we determine, for example, best possible value of $r_0$ so that
%$\|T_{f_r}\|_1\leq 1$ when $f\in \es$ and $0<r\leq r_0$; best possible constant $C(\alpha, r)$ for which
%$\|S_{f_r}\|_\alpha\le C(\alpha, r)$ holds for $f\in\es$, $\alpha >0$ and $0<r<1$.
%We obtain the largest possible $r\in (0,1)$ for which the pre-Schwarzian
%norm of the dilation $r^{-1}f(rz)$ is less than or equal to $1$,
%where $f$ is univalent in the unit disk with usual normalization.
%Similar results concerning Schwarzian derivative are also obtained.
\end{abstract}

\thanks{${}^{\mathbf{*}}$
T.~S.~was supported in part by JSPS Grant-in-Aid for
Scientific Research (B) 22340025.
}
\maketitle
\pagestyle{myheadings}
\markboth{S. Ponnusamy,  S.K. Sahoo and T. Sugawa}{(Pre-)Schwarzian
derivative and radius problems}

\section{Introduction}
Let $\A$ (resp.~$\M$)
denote the set of analytic (resp.~meromorphic) functions $f$ on the unit disk
$\D=\{z\in\C:|z|<1\}$ normalized so that $f(0)=0$ and $f'(0)=1.$
The set $\es$ of univalent functions in $\A$ has been intensively
studied by many authors.
The subclass of $\es$ consisting of convex functions
(i.e., functions mapping $\D$ univalently onto convex domains)
is denoted by $\K$, and the subclass of starlike functions
(i.e., functions mapping $\D$ univalently onto domains starlike with respect
to the origin) is denoted by $\es^\star.$
Thus, $\K\subset \es^\star \subset \es.$
Let $\mathcal F$ and $\mathcal G$ be two subclasses of $\A$.
If for every $f\in {\mathcal F}$, $r^{-1} f(rz)\in\mathcal{G}$ for
$0<r\leq r_0$, and
$r_0$ is the maximum value for which this holds, then we say that $r_0$ is
the $\mathcal{G}$-radius of $\mathcal{F}$. There are many results of this type
in the theory of univalent functions.
See \cite{Goodman:univ} for vast information in this direction.
For example, the $\K$-radius of $\mathcal{S},$
which is usually called the \textit{radius of convexity} of $\es,$
is $r_0^K=2-\sqrt{3}= 0.267...,$ (cf.~\cite[Theorem 2.13]{Duren:univ}).
On the other hand,
${\mathcal S}^\star$-radius in $\mathcal{S}$
(called the \textit{radius of starlikeness} of $\es$)
is $r_0^\star= \tanh (\pi/4) = 0.655...,$
(cf.~\cite[p.~98]{Duren:univ}). Krzy\.z \cite{Kr62B}
applied Loewner's method to show
that the radius of close-to-convexity of $\es$ is $r_0$, where
the equation defining $r_0$ is transcendental and successive approximations
yield that the lower and upper estimates for $r_0$ is $0.80$ and $0.81$, respectively.

It is sometimes important to give univalence criteria
for a non-constant meromorphic function $f$ on $\D$
in terms of its pre-Schwarzian or Schwarzian derivatives,
that are defined by
$$
T_f(z)=\frac{f''(z)}{f'(z)}\quad\text{and}\quad
S_f(z)=T_f^{\,\prime}(z)-\frac12T_f(z)^2,
%=\left(\frac{f''}{f'}\right)'(z)
%-\frac12\left(\frac{f''(z)}{f'(z)}\right)^2,
$$
respectively.
Note that $T_f$ (resp.~$S_f$) is analytic on $\D$ precisely when
$f$ is analytic (resp.~meromorphic) and locally univalent on $\D.$
In the theory of Teichm\"uller spaces, these quantities are
considered as elements of complex Banach spaces as follows.
For $\alpha\ge0,$ we define the norm
$$
\|\varphi\|_\alpha=\sup_{z\in\D}(1-|z|^2)^\alpha|\varphi(z)|
$$
for a meromorphic function $\varphi$ on $\D.$
Here, we define $\|\varphi\|_\alpha=+\infty$ whenever $\varphi$
has a pole in $\D.$
We denote by $\Bers_\alpha$ the complex Banach space consisting of
analytic functions $\varphi$ on $\D$ with $\|\varphi\|_\alpha<\infty.$
It is known that $f\in\A$ (resp.~$f\in\M$) is uniformly locally univalent
if and only if $T_f\in\Bers_1$ (resp.~$S_f\in\Bers_2$).
Here, a function $f$ is called {\it uniformly locally univalent}
if there exists a positive constant $\rho=\rho(f)$ such that
$f(z)$ is univalent on the disk $|(z-a)/(1-\bar az)|<\rho$
for every $a\in\D$
(see, for instance, \cite{Sug07ut}).

For $\alpha\ge0$ and $M>0,$ we set
$$
\B_\alpha(M)=\{f\in \A: \|T_f\|_\alpha\le M\}
$$
and
$$
\N_\alpha(M)=\{f\in \M: \|S_f\|_\alpha\le M\}.
$$

As a consequence of the area theorem (see \cite[p.~32]{Duren:univ}),
for $f\in\es$ we have the inequality
$$
|(1-|z|^2)T_f(z)-2\bar z|\le 4,
$$
and therefore
\begin{equation}\label{eq:ps}
|T_f(z)|\le\frac{2(2+|z|)}{1-|z|^2}.
\end{equation}
The last inequality leads to the implication
$\es\subset\B_1(6).$
Note here that the Koebe function $k(z)=z/(1-z)^2$ satisfies the relation
$(1-z^2)T_k(z)=2z+4,$ which shows the inequality \eqref{eq:ps} is sharp.
On the other hand,
Becker \cite{Becker72} showed the remarkable fact that
$\B_1(1)\subset \es.$
Sharpness of the constant $1$ is due to Becker and Pommerenke \cite{BP84}.

For the case of Schwarzian derivative, Nehari's result
\cite{Nehari49} is fundamental:
$\N_2(2)\subset \es\subset \N_2(6).$
(The latter inclusion relation is sometimes called the Kraus-Nehari theorem.)
The set $\N_2(2)$ is called the Nehari class and intensively studied
by Chuaqui, Osgood and Pommerenke \cite{COP96} (see also \cite{CP99}).
Note that the Koebe function $k(z)=z/(1-z)^2$ and its rotations are
not contained in $\N_2(2)$ since $\|S_k\|_2=6.$
It is also known that $\K\subset\N_2(2)$ and that the constant $2$ is sharp
(cf.~\cite[Lemma 1]{Rob69}).

There are variations of this type.
Nehari \cite{Nehari49} proved also that $\N_0(\pi^2/2)\subset\es$
and Pokornyi \cite{Pokornyi51} claimed that $N_1(4)\subset\es$
(see also \cite{Nehari54}).
The constants $\pi^2/2$ and $4$ are sharp (see \cite{Chuaqui95} and
\cite[\S 8.5 or p.264]{Duren:univ}, respectively).
For more refinements and background, see \cite{AA75}.

We begin the discussion with $f\in \es$ and its dilations
$$
f_r(z)=\frac{1}{r}f(rz), \quad 0\le r< 1.
$$
Here, $f_0$ is defined as the limit of $f_r;$ namely, $f_0(z)=z.$
Each function $f_r$, together with $f$, evidently belongs to $\es.$
Moreover, the relations
$$T_{f_r}(z)=rT_f(rz) ~\mbox{ and }~ S_{f_r}(z)=r^2S_f(rz),
$$
lead to the inequalities
$$\|T_{f_r}\|_\alpha \le \|T_f\|_\alpha~\mbox{ and }~\|S_{f_r}\|_\alpha \le \|S_f\|_\alpha
~\mbox{ for $0\le r<1$ and $\alpha\ge0.$}
$$
Since $T_{f_0}=S_{f_0}=0,$ one expects that
$\|T_{f_r}\|_\alpha\to0$ and $\|S_{f_r}\|_\alpha\to0$ as $r\to0.$
Thus, it may be interesting to estimate $\|T_{f_r}\|_\alpha$
and $\|S_{f_r}\|_\alpha$ and, moreover, for a given $M>0,$
to find the largest possible values $r_0$ of $r$ such that $f_r\in \B_\alpha(M)$
or $f_r\in \N_\alpha(M)$ for all $f\in\es.$
Note that $r_0$ is the $\B_\alpha(M)$-radius or the $\N_\alpha(M)$-radius
of $\es,$ respectively.

\section{Pre-Schwarzian derivative and radius property}

We first show the following result.

\begin{thm}\label{thm:ps1}
Let $0<r<1$ and $0<\alpha.$
Set $f_r(z)=r\inv f(rz)$ for an $f\in\es.$
Then $\|T_{f_r}\|_\alpha\le H(x_0),$ where
$$
H(x)=\frac{2r(1-x^2)^\alpha(2+rx)}{1-r^2x^2}
$$
and $x_0$ is the unique root of the polynomial
$$
Q(x)=r-4(\alpha-r^2)x-(1-r^2+2\alpha)rx^2+4(\alpha-1)r^2x^3+(2\alpha-1)r^3x^4
$$
in the interval $0<x<1.$
Moreover, equality holds when $f$ is the Koebe function.
\end{thm}
\begin{pf}
By \eqref{eq:ps}, we have the estimate
$$
(1-|z|^2)^\alpha |T_{f_r}(z)|
=(1-|z|^2)^\alpha r|T_{f}(rz)|
\le (1-|z|^2)^\alpha\frac{2r(2+r|z|)}{1-r^2|z|^2}
=H(|z|).
$$
Hence, the problem reduces to finding the supremum of $H(x)$ over $0<x<1.$
We now observe the formula
$$
\frac{H'(x)}{H(x)}=\frac{Q(x)}{(1-x^2)(2+rx)(1-r^2x^2)}.
$$
Since $Q(0)=r>0, Q(1)=-2\alpha(1-r^2)(2+r)<0,$ the intermediate value theorem
guarantees existence of a root of $Q(x)$ in the interval $0<x<1.$
Thus, it is enough to check uniqueness of the root in $0<x<1.$
To this end, we look at
$$
Q''(x)=2r\big[6r^2(2\alpha-1)x^2+12(\alpha-1)rx-(1-r^2+2\alpha)\big].
$$
When $\alpha\le 1/2,$ obviously $Q''(x)<0$ for $0<x<1,$ which implies
that $Q(x)$ is concave there. Hence, the root is unique.

We now suppose that $\alpha>1/2.$
Note first that $Q''(0)<0.$
Let $x_1$ be the (unique) positive root of the quadratic polynomial $Q''(x).$
When $x_1\ge1,$ $Q''(x)<0$ for $0<x<1.$
Thus $Q(x)$ is concave in $0<x<1$ which is enough to see the uniqueness.
When $x_1<1,$ $Q(x)$ is concave in $0<x<x_1$ and $Q(x)$ is convex in $x_1<x<1.$
If $Q(x_1)<0,$ then $Q(x)$ has a unique root in $0<x<x_1$ and has no root
in $x_1\le x<1.$
If $Q(x_1)\ge0,$ then $Q(x)$ has no root in $0<x<x_1$ and has a unique root
in $x_1\le x<1.$
At any event, the root is unique.
\end{pf}

\begin{rem}
When $\alpha=0,$ $Q(x)=r(1-x^2)(1+4rx+r^2x^2)\ge0.$
Thus,
$$H(x)<H(1)=\frac{2r(2+r)}{1-r^2}
$$
for $0\le x<1$ in this case.
Consequently, we have the sharp estimate
$$\|T_{f_r}\|_0\le \frac{2r(2+r)}{1-r^2}.
$$
\end{rem}

As an application of the last theorem, at least in principle,
we could find the $\B_\alpha(M)$-radius of $\es.$
As a simple example, we obtain the following.

\begin{thm}\label{thm:ps}
Let $r_0$ be the $\B_1(1)$-radius of $\es,$ that is,
the largest possible value so that $\|T_{f_r}\|_1\le 1$ whenever
$f\in\es$ and $0<r\le r_0.$
Then $r_0$ is the unique root of the equation
$$
17r^5-84r^4-61r^3-60r^2+277r-64=0
$$
in the interval $0<r<1$ and it is approximately $0.2489802.$
\end{thm}
\begin{pf}
It is enough to solve the system of equations
$$
H(r)=\frac{2r(1-x^2)(2+rx)}{1-r^2x^2}=1,
$$
and
$$
Q_1(x,r):=r^3x^4-(3-r^2)rx^2-4(1-r^2)x+r=0
$$
in the range $0<r<1$ and $0<x<1.$
Since the bound in Theorem \ref{thm:ps} is obviously increasing in $0<r<1,$
such a solution is unique in this range.
Let
$$
Q_2(x,r)=2r(1-x^2)(2+rx)-(1-r^2x^2).
$$
Then our task is to find a common zero of $Q_1$ and $Q_2.$
We divide $Q_1$ by $Q_2$ with respect to $x$ to obtain
$$
4Q_1(x,r)=(-2rx+4-r)Q_2(x,r)+Q_3(x,r),
$$
where
$$
Q_3(x,r)=(9r^3-8r^2+4r)x^2+(2r^3+16r^2-2r-16)x+4r^2-13r+4.
$$
Note that the set of common zeros of $Q_1$ and $Q_2$ is the same as
that of $Q_2$ and $Q_3.$
We repeat this procedure to have
$$
(9r^2-8r+4)^2Q_2(x,r)=\big(-2r(9r^2-8r-4)x+13r^3-12r^2+32r-48\big)
Q_3(x,r)+Q_4(x,r),
$$
where
$$
Q_4(x,r)=-8(1-r^2)\big((17r^4-50r^3+33r^2-56r+96)x+34r^3-55r^2-78r-22\big).
$$
By solving $Q_4(x,r)=0,$ we have the relation
$$
x=\frac{-34r^3+55r^2+78r+22}{17r^4-50r^3+33r^2-56r+96}.
$$
We substitute it into $Q_3(x,r):$
$$
Q_3(x,r)=-\frac{3(9r^2-8r+4)^2(17r^5-84r^4-61r^2+277r-64)}%
{(17r^4-50r^3+33r^2-56r+96)^2}.
$$
Thus we conclude that $r_0$ is a root of the equation in the assertion.
By the uniqueness of a solution to the system in the range $0<r<1, 0<x<1,$
we see that such a root is unique in the interval $0<r<1.$
Thus the proof is complete.
\end{pf}

%\vspace{6pt}

It is also important to observe what happens in Theorem \ref{thm:ps1} as $r\to1.$
Let
$$
P(\alpha)=\sup_{f\in\es}\|T_f\|_\alpha.
$$
Then we have the following result.

\begin{thm}
$$
P(\alpha)=\begin{cases}
+\infty
&\quad \text{if}\quad 0<\alpha<1, \\
6
&\quad \text{if}\quad \alpha=1, \\
\dfrac{(2\alpha+\sqrt{4\alpha^2-6\alpha+3})(3\alpha-3)^{\alpha-1}}
{(\alpha-1/2)(\alpha-3/2+\sqrt{4\alpha^2-6\alpha+3})^{\alpha-1}}
&\quad \text{if}\quad 1<\alpha .
\end{cases}
$$
\end{thm}
\begin{pf}
When $\alpha =1,$ the result is well-known and for $\alpha<1,$ the result is obvious.
We thus assume that $\alpha>1.$ Let
$$H(x)=2(1-x^2)^{\alpha-1}(2+x).
$$
Then, as in the proof of Theorem \ref{thm:ps1}, $P(\alpha)$ is given
as the maximum of $H(x)$ over $0<x<1.$
We now have
$$
\frac{H'(x)}{H(x)}=-\frac{(2\alpha-1)x^2+4(\alpha-1)x-1}{(1-x^2)(2+x)}
=-\frac{Q(x)}{(1-x^2)(2+x)}.
$$
Since $Q(0)=-1<0$ and $Q(1)=6(\alpha-1)>0,$ we have the unique root
$x_0$ of $Q(x)$ in $0<x<1$ and it is indeed given by
$$
x_0=\frac{\sqrt{4\alpha ^2-6\alpha+3}-2(\alpha-1)}{2\alpha-1}.
$$
Noting the relation
$$(2x_0+1)\left (2x_0+\frac{6\alpha-7}{2\alpha-1}\right )=3,
$$
we have
$$
1-x_0^2=\frac{2(\alpha-1)(2x_0+1)}{2\alpha-1}
=\frac{6(\alpha-1)}{2(2\alpha-1)x_0+6\alpha-7}
=\frac{6(\alpha-1)}{2\alpha-3+2\sqrt{4\alpha^2-6\alpha+3}}.
$$
Hence, $P(\alpha)=H(x_0)=2(1-x_0^2)^{\alpha-1}(2+x_0)$
has the form given in the assertion.
\end{pf}

%\bigskip

By definition, $P(\alpha)$ is non-increasing in $\alpha.$
Moreover, it is easy to see that $P(\alpha)\to 4$ as $\alpha\to+\infty.$
We also remark that the counterpart $\tilde P(\alpha)$ to the Schwarzian
derivative is very simple. Indeed, $\tilde P(\alpha)=+\infty$ if
$\alpha<2$ and $\tilde P(\alpha)=6$ otherwise.

\section{Schwarzian derivative and radius properties}

In this section we consider the Schwarzian derivative and its norm.
Our aim is to find the best possible constant $C(\alpha, r)$ for which
$\|S_{f_r}\|_\alpha\le C(\alpha, r)$ holds for $f\in\es.$
The result can be stated in the following form.

\begin{thm}\label{thm:s}
Let $0<r<1$ and $0<\alpha.$
Set $f_r(z)=r\inv f(rz)$ for an $f\in\es.$
Then
$$
\|S_{f_r}\|_\alpha\le\begin{cases}
%\displaystyle
6r^2\left(\dfrac{\alpha}{2r^2}\right)^\alpha
\left(\dfrac{1-\alpha/2}{1-r^2}\right)^{2-\alpha} &\quad
\text{if}\quad 0<\alpha<2r^2 \\
6r^2 &\quad \text{if}\quad 2r^2\le\alpha.
\end{cases}
$$
Equality holds when $f$ is the Koebe function.
\end{thm}
\begin{pf}
We recall the Kraus-Nehari theorem:
\begin{equation}\label{eq:KN}
(1-|z|^2)^2|S_f(z)|\le 6, \quad z\in \ID,
\end{equation}
for $f\in \es.$
It is a simple exercise to see that
$$(1-|z|^2)^\alpha|S_{f_r}(z)|=(1-|z|^2)^\alpha r^2|S_f(rz)|\le 6r^2K(|z|^2),
$$
where
$$
K(t)=\frac{(1-t)^\alpha}{(1-r^2t)^2}.
$$

Thus, it is sufficient to find the supremum of the function $K(t)$
over $0<t<1.$
We first look at the formula
$$
\frac{K'(t)}{K(t)}
=\frac{(\alpha-2)r^2t+2r^2-\alpha}{(1-t)(1-r^2t)}
=-\frac{(\alpha-2r^2)(1-r^2t)+2r^2(1-r^2)}{(1-t)(1-r^2t)}.
$$
When $\alpha\ge 2r^2,$ obviously $K'(t)\le0$ in $0<t<1.$
Therefore, in this case, $K(t)$ is non-increasing in $t$ and its supremum
is $K(0)=1.$
On the other hand, when $0<\alpha<2r^2,$ the function $K(t)$
takes its maximum at
$$
t_0
=\frac{2r^2-\alpha}{(2-\alpha)r^2}
=\frac{2r^2-\alpha}{2r^2-\alpha r^2}.
$$
Since
$$
K(t_0)
=\left(\frac{\alpha(1-r^2)}{r^2(2-\alpha)}\right)^\alpha
\left(\frac{2-\alpha}{2(1-r^2)}\right)^2,
$$
the conclusion follows.
\end{pf}

%\vspace{6pt}

The above theorem determines the value of
$$
C_\alpha(r)=\sup_{f\in\es}\|S_{f_{r}}\|_\alpha.
$$
When $\alpha=0,$ the above computation tells us that $K(t)$ is increasing
in $0<t<1.$
Therefore, $C_0(r)=6r^2K(1)=6r^2(1-r^2)^{-2}.$
We thus summarize the conclusions:
$$
C_\alpha(r)=\begin{cases}
\dfrac{6r^2}{(1-r^2)^2} &\quad
\text{if}\quad \alpha=0 \\
6r^2\left(\dfrac{\alpha}{2r^2}\right)^\alpha
\left(\dfrac{1-\alpha/2}{1-r^2}\right)^{2-\alpha} &\quad
\text{if}\quad 0<\alpha<2r^2 \\
6r^2 &\quad \text{if}\quad 2r^2\le\alpha.
\end{cases}
$$
In particular, we observe that $C_\alpha(r)$
is strictly increasing in $0<r<1.$
In particular, we have the following.

\begin{cor}
Let $0<M\le 3\alpha.$
Then the $\N_\alpha(M)$-radius of $\es$ is $\sqrt{M/6}.$
\end{cor}
\begin{pf}
Let $r_0=\sqrt{M/6}.$
Then $2r_0^2=M/3\le\alpha$ and thus Theorem \ref{thm:s} implies
$C_\alpha(r_0)=6r_0^2=M.$
\end{pf}

%\vspace{8pt}

For instance, we easily have the following:
\begin{enumerate}
\item The $\N_2(2)$-radius of $\es$ is $1/\sqrt{3}=0.5773\dots.$
\item  The $\N_0(M)$-radius of $\es$ is obtained by solving the equation
$$C_0(r)=\frac{6r^2}{(1-r^2)^2}=M
$$
which gives
$$
r=\sqrt{\frac{M+3-\sqrt{9+6M}}{M}}.
$$
\item The $\N_0(\pi^2/2)$-radius of $\es$ is
$\sqrt{1+(6/\pi^2)-(2\sqrt{9+3\pi^2}/\pi^2)}
=0.5905\dots.$
\item Finally, we find $\N_1(4)$-radius of $\es.$
In this case, we solve the equation
$$C_1(r)=\frac{3}{2(1-r^2)}=4
$$
to obtain $r=\sqrt{5/8}=0.7905\dots.$
\end{enumerate}

\def\cprime{$'$} \def\cprime{$'$} \def\cprime{$'$}
\providecommand{\bysame}{\leavevmode\hbox to3em{\hrulefill}\thinspace}
\providecommand{\MR}{\relax\ifhmode\unskip\space\fi MR }
\providecommand{\MRhref}[2]{%
  \href{http://www.ams.org/mathscinet-getitem?mr=#1}{#2}
}
\providecommand{\href}[2]{#2}

%\bibliography{papers}
\end{document}